\pdfoutput=1
\documentclass[toc]{mathprint}
\usepackage{rutarmacros}
\usepackage{project-macros}
\title{Tangents of invariant sets}
\author
  {Antti Käenmäki}
  {Alfréd Rényi Institute of Mathematics,
   Hungarian Academy of Sciences,
   Budapest,
   Hungary}
  {kaenmaki@renyi.hu}

\author
  {Alex Rutar}
  {Mattilanniemi (MaD),
   40100 University of Jyväskylä,
   Finland}
  {alex@rutar.org}

\begin{document}

\begin{abstract}
    We study the fine scaling properties of sets satisfying various weak forms of invariance.
    For general attractors of possibly overlapping bi-Lipschitz iterated function systems, we establish that the Assouad dimension is given by the Hausdorff dimension of a tangent at some point in the attractor.
    Under the additional assumption of self-conformality, we moreover prove that this property holds for a subset of full Hausdorff dimension.
\end{abstract}

\section{Introduction}
One of the most fundamental concepts at the intersection of analysis and geometry is the notion of a \emph{tangent}.
For sets exhibiting a high degree of local regularity---such as manifolds, or rectifiable sets---at almost every point in the set and at all sufficiently high resolutions, the set looks essentially linear.
Moreover, the concept of a tangent is particularly relevant in the study of a much broader class of sets: those equipped with some form of dynamical invariance.
This relationship originates in the pioneering work of Furstenberg, where one associates to a set a certain dynamical system of ``zooming in''.
Especially in the past two decades, the study of tangent measures has played an important role in the resolution of a number of long-standing problems concerning sets which look essentially the same at all small scales; see, for example, \cite{zbl:1409.11054,zbl:1251.28008,zbl:1430.11106,zbl:1426.11079,zbl:1318.28026}.

In contrast, (weak) tangents also play an important role in the geometry of metric spaces.
One of the main dimensional quantities in the context of embeddability properties of metric spaces is the Assouad dimension, first introduced in \cite{ass1977}.
It turns out that the Assouad dimension, which bounds the worst-case scaling at all locations and all small scales, is precisely the maximal Hausdorff dimension of weak tangents, i.e.~sets which are given as a limit of small pieces of enlarged copies of the original set; see \cite{doi:10.1093/imrn/rnw336}.
We refer the reader to the books \cite{fra2020,rob2011,mt2010} for more background and context on the importance of Assouad dimension in a variety of diverse applications.

In this document, we study the interrelated concepts of tangents and Assouad dimension, with an emphasis on sets with a weak form of dynamical invariance.
Our motivating examples include attractors of iterated function systems where the maps are affinities (or even more generally bi-Lipschitz contractions); or the maps are conformal and there are substantial overlaps.
In both of these situations, the sets exhibit a large amount of local inhomogeneity.

\subsection{Weak tangents, tangents, and pointwise Assouad dimension}
Throughout, we will work in $\R^d$ for some $d\in\N$, though many of our results hold in the broader context of bounded doubling metric spaces.
We let $B(x,r)$ denote the closed ball with centre $x$ and radius $r$.

Now, fix a compact set $K\subset\R^d$.
We say that a compact set $F\subset B(0,1)$ is a \defn{weak tangent} of $K\subset\R^d$ if there exists a sequence of similarity maps $(T_k)_{k=1}^\infty$ with similarity ratios $\lambda_k$ diverging to infinity such that $0\in T_k(K)$ and
\begin{equation*}
    F=\lim_{k\to\infty}T_k(K)\cap B(0,1)
\end{equation*}
with respect to the Hausdorff metric on compact subsets of $B(0,1)$.
We denote the set of weak tangents of $K$ by $\Tan(K)$.
More strongly, we say that $F$ is a \defn{tangent of $K$ at $x$} if $F$ is a weak tangent and the similarity maps $T_k$ are homotheties which map $x$ to $0$; i.e.\ $T_k(y)= \lambda_k(y-x)$.
We denote the set of tangents of $K$ at $x$ by $\Tan(K,x)$.
We refer the reader to \cref{ss:tan-def} for precise definitions.

Closely related to the notion of a weak tangent is the \defn{Assouad dimension} of $K$, which is the dimensional quantity
\begin{align*}
    \dimA K=\inf\Bigl\{s:\exists C>0\,&\forall 0<r\leq R<1\,\forall x\in K\\*
                                      &N_r(B(x,R)\cap K)\leq C\Bigl(\frac{R}{r}\Bigr)^s\Bigr\}.
\end{align*}
Here, for a general bounded set $F$, $N_r(F)$ is the smallest number of closed balls with radius $r$ required to cover $F$.
It always holds that $\dimH K\leq\dimuB K\leq\dimA K$, where $\dimH K$ and $\dimuB K$ denote the Hausdorff and upper box dimensions respectively.
In some sense, the Assouad dimension is the largest reasonable notion of dimension which can be defined using covers.
Continuing the analogy with tangents, we also introduce a localized version of the Assouad dimension which we call the \defn{pointwise Assouad dimension}.
Given $x\in K$, we set
\begin{align*}
    \dimA(K,x)=\inf\Bigl\{s:\exists C>0\,&\exists\rho>0\,\forall 0<r\leq R<\rho\\*
                                         &N_r(B(x,R)\cap K)\leq C\Bigl(\frac{R}{r}\Bigr)^s\Bigr\}.
\end{align*}
The choice of $\rho>0$ in the definition of $\dimA(K,x)$ ensures a sensible form of bi-Lipschitz invariance: if $f\colon K\to K'$ is bi-Lipschitz, then $\dimA(K,x)=\dimA(f(K),f(x))$.
It is immediate from the definition that
\begin{equation*}
    \dimA(K,x)\leq\dimA K.
\end{equation*}
Moreover, if for instance $K$ is Ahlfors--David regular, then $\dimA(K,x)=\dimA K$ for all $x\in K$.
We note here that an analogous notion of pointwise Assouad dimension for measures was introduced recently in \cite{zbl:1527.28004}.

An important observation which is essentially due to Furstenberg \cite{zbl:0208.32203,zbl:1154.37322}, but was observed explicitly in \cite{doi:10.1093/imrn/rnw336}, is that the Assouad dimension is characterized by weak tangents:
\begin{equation*}
    \dimA K=\max\{\dimH F:F\in\Tan(K)\}.
\end{equation*}
Motivated by this relationship, our primary goal in this document is to address the following questions:
\begin{itemize}[nl]
    \item Does it hold that $\dimA(K,x)=\max\{\dimH F:F\in\Tan(K,x)\}$?
    \item Is there necessarily an $x_0\in K$ so that $\dimA K=\dimH F$ for some $F\in\Tan(K,x_0)$?
        If not, is there an $x_0\in K$ so that $\dimA K=\dimA(K,x_0)$?
    \item What is the structure of the level set of pointwise Assouad dimension $\{x\in K:\dimA(K,x)=\alpha\}$ for some $\alpha\geq 0$?
\end{itemize}
In the following section, we discuss our main results and provide some preliminary answers which indicate that answers to these questions are, in general, quite subtle.

\subsection{Main result and outline of paper}
We begin by stating some easy properties of the pointwise Assouad dimension for general compact sets $K\subset\R^d$.
First, some standard measurability results are stated in \cref{p:meas}.
Next, by \cref{p:pointwise-tangent},
\begin{equation*}
    \sup\{\dimuB F:F\in\Tan(K,x)\} \leq \dimA(K,x) \leq \dimA K.
\end{equation*}
Unfortunately, in general one cannot hope for either inequality to be in equality: an example in \cite{zbl:1321.54059} already has the property that $K\subset\R$ such that $\dimA K =1$ but $\dimA(K,x)=0$ for all $x\in K$ (see \cref{ex:non-tan-attain} for more detail); and moreover, in \cref{ex:non-point-attain}, we construct a compact set $K\subset\R$ with a point $x\in K$ so that $\dimA(K,x)=1$ but each $F\in\Tan(K,x)$ consists of at most two points.

In light of the above, for general sets, the following result is essentially optimal.
\begin{itheorem}\label{it:gen-tangents}
    Let $K\subset\R^d$.
    Then:
    \begin{enumerate}[nl,r]
        \item\label{im:1} If $K$ is compact, there is an $x\in K$ such that $\dimA(K,x)\geq\dimuB K$.
        \item\label{im:2} If $K$ is analytic, for any $s$ such that $\mathcal{H}^s(K)>0$, there is a compact set $E\subset K$ with $\mathcal{H}^s(E)>0$ so that for each $x\in E$, there is a tangent $F\in\Tan(\overline{K},x)$ with $\mathcal{H}^{s}_\infty(F)\geq 1$.
    \end{enumerate}
\end{itheorem}
The proof of \cref{im:1} is a direct consequence of the definitions.
Moreover, \cref{im:2} is almost a consequence of the classical density theorem for Hausdorff measure.
However, the main difficulty is that Hausdorff measure and the Hausdorff metric are in general not compatible.
To work around this, we require a simple but key observation: the density theorems also hold for Hausdorff content, and the Hausdorff content map $K\mapsto\mathcal{H}^s_\infty(K)$ is upper semicontinuous with respect to the Hausdorff metric on the space of compact sets.

Actually, \cref{it:gen-tangents} has a useful consequence for general sets, which will also play a key role in the proof of \cref{it:self-embed} below.
\begin{icorollary}\label{ic:full-tan}
    Let $K\subseteq\R^d$ be a non-empty compact set with $\alpha=\dimA K$.
    Then there is an $F\in\Tan(K)$ such that $\mathcal{H}^\alpha_\infty(F)\geq 1$.
\end{icorollary}

While we cannot hope to improve \cref{it:gen-tangents} in general, many commonly studied families of ``fractal'' sets have a form of dynamical invariance, which is far from the case for general sets.
As a result, it is of interest to determine general conditions under which the Assouad dimension is actually attained as the pointwise Assouad dimension at some point.
To this end, we make the following definition.
\begin{definition}\label{d:self-embed}
    We say that a compact set $K$ is \defn{self-embeddable} if for each $z\in K$ and $0<r\leq\diam K$, there is a constant $a=a(z,r)>0$ and a function $f\colon K\to B(z,r)\cap K$ so that
    \begin{equation}\label{e:dist-bound}
        a r|x-y|\leq |f(x)-f(y)|\leq a^{-1}r|x-y|.
    \end{equation}
    for all $x,y\in K$.
    We say that $K$ is \defn{uniformly self-embeddable} if the constant $a(z,r)$ can be chosen independently of $z$ and $r$.
\end{definition}
The class of self-embeddable sets is very broad and includes, for example, attractors of every possibly overlapping iterated function system $\{f_i\}_{i\in\mathcal{I}}$, where $\mathcal{I}$ is a finite index set and $f_i$ is a strictly contracting bi-Lipschitz map from $\R^d$ to $\R^d$.

The class of \emph{uniformly} self-embeddable sets includes the attractors of finite overlapping self-conformal iterated function systems.
It is perhaps useful to compare uniform self-embeddability with quasi self-similarity, as introduced by Falconer \cite{zbl:0683.58034}.
Our assumption is somewhat stronger since we also require the upper bound to hold in \cref{e:dist-bound}.
This assumption is critical to our work since, in general, maps satisfying only the lower bound can decrease Assouad dimension.
We also note that uniform self-embeddability is the primary assumption in \cite[Theorem~2.1]{zbl:1441.28003}.

Within this general class of sets, we establish the following result which guarantees the existence of at least one large tangent under self-embeddability, and an abundance of tangents under uniform self-embeddability.
\begin{itheorem}\label{it:self-embed}
    Let $K\subset\R^d$ be compact and self-embeddable.
    Then:
    \begin{enumerate}[nl,r]
        \item $\dimuB K\leq\dimA(K,x)$ for all $x\in K$.
        \item\label{im:se} There is a dense set $P\subset K$ such that for each $x\in P$, there is a $F\in\Tan(K,x)$ so that $\mathcal{H}^{\dimA K}(F)\geq 2^{-\dimA K}$.
            In particular,
            \begin{equation}\label{e:pack-full}
                \dimP\{x\in K:\dimA(K,x)=\dimA K\}=\dimP K.
            \end{equation}
    \end{enumerate}
    If $K$ is uniformly self-embeddable, then there is a constant $c>0$ so that
    \begin{equation}\label{e:uniform-bd}
        \dimH\{x\in K:\exists F\in\Tan(K,x)\text{ with }\mathcal{H}^{\dimA K}(F)\geq c\}=\dimH K.
    \end{equation}
\end{itheorem}
\cref{it:self-embed} can be obtained by combining \cref{t:pointwise-tan}, \cref{p:tan-lower}, and \cref{t:uniformly-self-embeddable}.
As a special case of the result for uniformly self-embeddable sets, suppose $K$ is the attractor of a finite self-similar IFS in the real line with Hausdorff dimension $s<1$.
In this case there is a dichotomy: either $\mathcal{H}^s(K)>0$, in which case $K$ is Ahlfors--David regular, or $\dimA K=1$.
In particular, \cref{e:uniform-bd} cannot be improved in general to give a set with positive Hausdorff $s$-measure.

Beyond being of general interest, we believe this result will be a useful technical tool in the study of Assouad dimension for general attractors of bi-Lipschitz invariant sets.
For instance, a common technique in studying attractors of iterated function systems is to relate the underlying geometry to symbolic properties associated with the coding space.
Upper bounding the Hausdorff dimension of tangents is \emph{a priori} easier since one may fix in advance a coding for the point.
This is the situation, for example, in \cite[Theorem~5.2]{zbl:1455.28005}.

\subsection{Some variants for future work}
Firstly, a natural question is to what extend the results in \cref{it:self-embed} can be improved.
More precisely, under the general self-embeddable assumption, can one replace packing dimension with Hausdorff dimension in \cref{e:pack-full}?
In subsequent work by the authors, as a consequence of a more detailed study of the structure of tangents of a certain family of self-affine carpets, we will give an example of a self-embeddable set for which
\begin{equation*}
    \dimH\{x\in K:\dimA(K,x)=\dimA K\}<\dimH K.
\end{equation*}
Of course, one might still wonder in general if there are natural conditions which are weaker than self-conformality under which \cref{e:uniform-bd} holds.

Secondly, we may also define a more general variant of the pointwise Assouad dimension.
Let $\phi\colon(0,1)\to(0,1)$ be a fixed function.
We then define the \defn{pointwise $\phi$-Assouad dimension}, given by
\begin{align*}
    \dimAs\phi(K,x)=\inf\Bigl\{s:\exists C>0\,&\forall 0<r<1\\*
                                              &N_{r^{1+\phi(r)}}\bigl(B(x,r)\cap K\bigr)\leq Cr^{-\phi(r) s}\Bigr\}.
\end{align*}
It is a straightforward to see that
\begin{equation*}
    \dimAs\phi(K,x)=\limsup_{r\to 0}\frac{\log N_{r^{1+\phi(r)}}\bigl(B(x,r)\cap K\bigr)}{\phi(r)\log(1/r)}.
\end{equation*}
The $\phi$-Assouad dimensions are an example of \emph{dimension interpolation} \cite{zbl:1462.28007} and have been studied in detail in \cite{ghm2021,arxiv:2308.12975}.
In the specific case that $\phi(R)=\frac{1}{\theta}-1$ for some $\theta\in(0,1)$, this corresponds precisely to the Assouad spectrum \cite{fy2018b} which (abusing notation) we may denote by $\dimAs\theta(K,x)$.
In general, we expect the properties of the pointwise Assouad spectrum to be substantially different than the properties of the pointwise Assouad dimension.

Thirdly, one might also consider the dual notion of the \defn{pointwise lower dimension}, defined for $x\in K$ by
\begin{align*}
    \dimL(K,x)=\sup\Bigl\{s:\exists C>0\,&\exists\rho>0\,\forall 0<r\leq R<\rho\\*
                                         &N_r(B(x,R)\cap K)\geq C\Bigl(\frac{R}{r}\Bigr)^s\Bigr\}.
\end{align*}
It is established in \cite{zbl:1428.28013} that the lower dimension may be analogously characterized as the minimum of Hausdorff dimensions of weak tangents.
Therefore, a natural question is to ask if similar results hold for the pointwise lower dimension as well.
However, the proofs we have given for \cref{it:self-embed} do not immediately translate to the case of the lower dimension since overlaps may increase dimension.

Finally, we note that an analogous notion for the pointwise Assouad dimension of measures was recently introduced in \cite{zbl:1527.28004}.
It would be interesting to investigate the relationship between these two notions of pointwise dimension.

\subsection{Notation}
Throughout, we work in $\R^d$ equipped with the usual Euclidean metric.
Write $\R_+=(0,\infty)$.
Given functions $f$ and $g$, we say that $f\lesssim g$ if there is a constant $C>0$ so that $f(x)\leq C g(x)$ for all $x$ in the domain of $f$ and $g$.
We write $f\approx g$ if $f\lesssim g$ and $g\lesssim f$.

\section{Tangents and pointwise Assouad dimension} \label{sec:pointwise}
\subsection{Tangents and weak tangents}\label{ss:tan-def}
To begin this section, we precisely define the notions of tangent and weak tangent, and establish the fundamental relationship between the dimensions of tangents and the pointwise Assouad dimension.

Given a set $E\subset\R^d$ and $\delta>0$, we denote the \emph{open $\delta$-neighbourhood} of $E$ by
\begin{equation*}
    E^{(\delta)}=\{x\in\R^d:\exists y\in E\text{ such that }|x-y|<\delta\}.
\end{equation*}
Now given a non-empty subset $X\subset\R^d$, we let $\mathcal{K}(X)$ denote the set of non-empty compact subsets of $X$ equipped with the \defn{Hausdorff metric}
\begin{equation*}
    d_{\mathcal{H}}(K_1,K_2)=\max\{p_{\mathcal{H}}(K_1;K_2),p_{\mathcal{H}}(K_2;K_1)\}
\end{equation*}
where
\begin{equation*}
    p_{\mathcal{H}}(K_1;K_2)=\inf\{\delta>0:K_1\subset K_2^{(\delta)}\}.
\end{equation*}
If $X$ is compact, then $(\mathcal{K}(X),d_{\mathcal{H}})$ is a compact metric space itself.
We also write
\begin{equation*}
    \dist(E_1,E_2)=\inf\{|x-y|:x\in E_1, y\in E_2\}
\end{equation*}
for non-empty sets $E_1,E_2\subset\R^d$.

We say that a set $F\in\mathcal{K}(B(0,1))$ is a \defn{weak tangent} of $K\subset\R^d$ if there exists a sequence of similarity maps $(T_k)_{k=1}^\infty$ with $0\in T_k(K)$ and similarity ratios $\lambda_k$ diverging to infinity such that
\begin{equation*}
    F=\lim_{k\to\infty}T_k(K)\cap B(0,1)
\end{equation*}
in $\mathcal{K}(B(0,1))$.
We denote the set of weak tangents of $K$ by $\Tan(K)$.
A key feature of the Assouad dimension is that it is characterized by Hausdorff dimensions of weak tangents.
This result is explicitly stated in \cite[Proposition~5.7]{doi:10.1093/imrn/rnw336}.
We refer the reader to \cite[Section~5.1]{fra2020} for more discussion on the context and history of this result.
\begin{proposition}[\cite{zbl:1154.37322,doi:10.1093/imrn/rnw336}]\label{p:tangent-assouad}
    We have
    \begin{equation*}
        \alpha\coloneqq\dimA K=\max_{F\in\Tan(K)}\dimH F.
    \end{equation*}
    Moreover, the maximizing weak tangent $F$ can be chosen so that $\mathcal{H}^{\alpha}(F)>0$.
\end{proposition}
In a similar flavour, we say that $F$ is a \defn{tangent of $K$ at $x\in K$} if there exists a sequence of similarity ratios $(\lambda_k)_{k=1}^\infty$ diverging to infinity such that
\begin{equation*}
    F=\lim_{k\to\infty}\lambda_k(K-x)\cap B(0,1)
\end{equation*}
in $\mathcal{K}(B(0,1))$.
We denote the set of tangents of $K$ at $x$ by $\Tan(K,x)$.

Of course, $\Tan(K,x)\subset\Tan(K)$.
Unlike in the case for weak tangents, we require the similarities in the construction of the tangent to in fact be homotheties.
This choice is natural since, for example, a function $f\colon\R\to\R$ is differentiable at $x$ if and only if the set of tangents of the graph of $f$ at $(x,f(x))$ is the singleton $\{B(0,1)\cap\ell\}$ for some non-vertical line $\ell$ passing through the origin.
In practice, compactness of the group of orthogonal transformations in $\R^d$ means this restriction will not cause any technical difficulties.

We observe that upper box dimensions of tangents provide a lower bound for the pointwise Assouad dimension.
\begin{proposition}\label{p:pointwise-tangent}
    For any compact set $K\subset\R^d$ and $x\in K$, $\dimA(K,x)\geq\dimuB F$ for any $F\in\Tan(K,x)$.
\end{proposition}
\begin{proof}
    Let $\alpha>\dimA(K,x)$ and suppose $F\in\Tan(K,x)$: we will show that $\dimuB F\leq\alpha$.
    First, get $C>0$ such that for each $0<r\leq R<1$,
    \begin{equation*}
        N_r(B(x,R)\cap K)\leq C\Bigl(\frac{R}{r}\Bigr)^\alpha.
    \end{equation*}
    Let $\delta>0$ be arbitrary, and get a similarity $T$ with similarity ratio $\lambda$ such that $T(x)=0$ and
    \begin{equation*}
        d_{\mathcal{H}}(T(K)\cap B(0,1),F)\leq\delta.
    \end{equation*}
    Then there is a uniform constant $M>0$ so that
    \begin{equation*}
        M\cdot N_\delta(F)\leq N_\delta(T(K)\cap B(0,1))=N_{\delta\lambda}(K\cap B(x,\lambda))\leq C\Bigl(\frac{\lambda}{\delta\lambda}\Bigr)^\alpha=C\delta^{-\alpha}.
    \end{equation*}
    In other words, $\dimuB F\leq\alpha$.
\end{proof}
One should not expect equality to hold in general: in \cref{ex:non-point-attain}, we construct an example of a compact set $K\subset\R$ and a point $x\in K$ so that $\dimA(K,x)=1$ but every $F\in\Tan(K,x)$ consists of at most 2 points.

\subsection{Level sets and measurability}\label{s:multi}
We now make some observations concerning the multifractal properties of the function $x\mapsto\dimA(K,x)$.
In particular, we are interested in the following quantities:
\begin{equation*}
    \mathcal{A}(K,\alpha)=\{x\in K:\dimA(K,x)=\alpha\}\qquad\text{and}\qquad\varphi(\alpha)=\dimH\mathcal{A}(K,\alpha).
\end{equation*}
We use the convention that $\dimH\varnothing=-\infty$.
Observe that $\varphi$ is a bi-Lipschitz invariant.

Let $\mathcal{K}(\R^d)$ denote the family of compact subsets of $\R^d$, equipped with the Hausdorff distance $d_{\mathcal{H}}$.
We recall that $B(x,r)$ denotes the closed ball at $x$ with radius $r$, and we let $B^\circ(x,r)$ denote the open ball at $x$ with radius $r$.
Given a compact set $K\subset\R^d$, we let $N^\circ_r(K)$ denote the minimal number of open sets with diameter $r$ required to cover $K$, and $N_r^{\mathrm{pack}}(K)$ denote the size of a maximal centred packing of $K$ by closed balls with radius $r$.
Then, for $0<r_1\leq r_2$, we write
\begin{align*}
    \mathcal{N}^\circ_{r_1,r_2}(K,x)&=N^\circ_{r_1}(B(x,r_2)\cap K)\\*
    \mathcal{N}_{r_1,r_2}(K,x)&=N_{r_1}^{\mathrm{pack}}(B^\circ(x,r_2)\cap K)
\end{align*}
The following lemma is standard.
\begin{lemma}
    Fix $0<r_1\leq r_2$.
    Then:
    \begin{enumerate}[nl,r]
        \item $\mathcal{N}^\circ_{r_1,r_2}\colon\mathcal{K}(\R^d)\times\R^d\to[0,d]$ is lower semicontinuous.
        \item $\mathcal{N}_{r_1,r_2}\colon\mathcal{K}(\R^d)\times\R^d\to[0,d]$ is upper semicontinuous.
    \end{enumerate}
\end{lemma}
We can use this lemma to establish the following fundamental measurability results.
\begin{proposition}\label{p:meas}
    The following measurability properties hold:
    \begin{enumerate}[nl,r]
        \item\label{im:large-gdelta} For a fixed compact set $K$ and $t\geq 0$, the set $\{x:\dimA(K,x)\geq s\}$ is a $G_\delta$ set.
        \item The function $(K,x)\mapsto\dimA(K,x)$ is Baire class 2.
        \item $\mathcal{A}(K,\alpha)$ is Borel for any compact set $K$.
    \end{enumerate}
\end{proposition}
\begin{proof}
    Since $\R^d$ is doubling,
    \begin{align*}
        \dimA(K,x)=\inf\Bigl\{s:\exists C>0\,&\exists M\in\N\,\forall M\leq k\leq n\\*
                                                                 &\mathcal{N}_{2^{-n},2^{-k}}(K,x)\leq C2^{(n-k)s}\Bigr\}.
    \end{align*}
    Equivalently, we may use $\mathcal{N}^\circ_{r_1,r_2}$ in place of $\mathcal{N}_{r_1,r_2}$.
    In particular,
    \begin{equation*}
        \{(K,x):\dimA(K,x)>s\}=\bigcap_{C=1}^\infty\bigcap_{M=1}^\infty\bigcup_{k=M}^\infty\bigcup_{n=k}^\infty(\mathcal{N}^\circ_{2^{-n},2^{-k}})^{-1}(C2^{(n-k)s},\infty)
    \end{equation*}
    is a $G_\delta$ set.
    Thus $\{x:\dimA(K,x)>s\}$ is also a $G_\delta$ set, so
    \begin{equation*}
        \{x:\dimA(K,x)\geq t\}=\bigcap_{n=1}^\infty\{x:\dimA(K,x)>t-1/n\}
    \end{equation*}
    is also a $G_\delta$ set, as claimed in \cref{im:large-gdelta}.

    Moreover,
    \begin{equation*}
        \{(K,x):\dimA(K,x)<t\}=\bigcup_{C\in\Q\cap(0,\infty)}\bigcup_{M=1}^\infty\bigcap_{k=M}^\infty\bigcap_{n=k}^\infty(\mathcal{N}_{2^{-n},2^{-k}})^{-1}(-\infty,C2^{(n-k)t}).
    \end{equation*}
    Thus $\{(K,x):\dimA(K,x)\in(s,t)\}$ is a $G_{\delta\sigma}$-set, i.e.~it is a countable union of sets expressible as a countable intersection of open sets, so $\dimA$ is Baire class 2.

    Of course, the same argument also show that $x\mapsto\dimA(K,x)$ is Baire class 2 for a fixed compact set $K$, so that $\mathcal{A}(K,\alpha)$ is $G_{\delta\sigma}$ and, in particular, Borel.
\end{proof}

\subsection{Tangents and pointwise dimensions of general sets}\label{ss:general-sets}
We now establish some general results on the existence of tangents for general sets.
These results will also play an important technical role in the following sections: for many of our applications, it is not enough to have positive Hausdorff $\alpha$-measure for $\alpha=\dimA K$, since in general Hausdorff $\alpha$-measure does not interact well with the Hausdorff metric on $\mathcal{K}\bigl(B(0,1)\bigr)$.

Recall that the \defn{Hausdorff $\alpha$-content} of a set $E$ is given by
\begin{equation*}
    \mathcal{H}^\alpha_\infty(E)=\inf\left\{\sum_{i=1}^\infty(\diam U_i)^\alpha:E\subset\bigcup_{i=1}^\infty U_i, U_i\text{ open}\right\}.
\end{equation*}
Of course, $\mathcal{H}^\alpha_\infty(E)\leq\mathcal{H}^\alpha(E)$ and $\mathcal{H}^\alpha_\infty(E)=0$ if and only if $\mathcal{H}^\alpha(E)=0$.
We recall (see, e.g.~\cite[Theorem~2.1]{zbl:0885.28005}) that $\mathcal{H}^\alpha_\infty$ is upper semicontinuous on $\mathcal{K}(B(0,1))$.
Moreover, if $0<\mathcal{H}^\alpha(E)<\infty$, then the density theorem for Hausdorff content implies that $\mathcal{H}^\alpha$-almost every $x\in E$ has a tangent with uniformly large Hausdorff $\alpha$-content.
We use these ideas in the following proofs.

We begin with a straightforward preliminary lemma which is proven, for example, in \cite[Lemma~3.11]{zbl:1342.28016}.
\begin{lemma}\label{l:tan-compose}
    Let $K\subset\R^d$ be compact.
    Then $\Tan(\Tan(K))\subset\Tan(K)$.
\end{lemma}
\begin{proof}
    First suppose $E\in\Tan(K)$ and $F\in\Tan(E)$.
    Write $E=\lim_{n\to\infty}T_n(K)\cap B(0,1)$ and $F=\lim_{n\to\infty}S_n(E)\cap B(0,1)$ for some sequences of similarities $(T_n)$ and $(S_n)$ with similarity ratios diverging to infinity.
    For each $\epsilon>0$, let $N$ be sufficiently large so that
    \begin{equation*}
        d_{\mathcal{H}}(S_N(E)\cap B(0,1),F)\leq\frac{\epsilon}{2}.
    \end{equation*}
    Suppose $S_N$ has similarity ratio $\lambda_N$, and let $M$ be sufficiently large so that
    \begin{equation*}
        d_{\mathcal{H}}(T_M(K)\cap B(0,1),E)\leq\frac{\epsilon}{2\lambda_N}.
    \end{equation*}
    It follows that
    \begin{equation*}
        d_{\mathcal{H}}(S_N\circ T_M(K)\cap B(0,1),F)\leq\epsilon.
    \end{equation*}
    But $\epsilon>0$ was arbitrary, as required.
\end{proof}

Now, given a set with positive and finite Hausdorff measure, we can always find a tangent with large Hausdorff content.
\begin{lemma}\label{l:large-tan}
    Let $K\subseteq\R^d$ be a compact set with $0<\mathcal{H}^\alpha(K)<\infty$.
    Then for $\mathcal{H}^\alpha$-almost every $x\in K$, there is an $F\in\Tan(K,x)$ such that $\mathcal{H}^\alpha_\infty(F)\geq 1$.
\end{lemma}
\begin{proof}
    By the same proof as \cite[Theorem~6.2]{zbl:0819.28004}, for $\mathcal{H}^\alpha$-almost every $x\in K$, there is a sequence of scales $(r_n)_{n=1}^\infty$ converging to zero such that
    \begin{equation*}
        1\leq\lim_{n\to\infty}r_n^{-\alpha}\mathcal{H}^\alpha_\infty\bigl(B(x,r_n)\cap K\bigr).
    \end{equation*}
    Then
    \begin{align*}
        \mathcal{H}^\alpha_\infty\bigl(r_n^{-1}(K-x)\cap B(0,1)\bigr)&= r_n^{-\alpha}\mathcal{H}^\alpha_\infty\bigl(B(x,r_n)\cap K\bigr)\fto{n\to\infty}1.
    \end{align*}
    But Hausdorff $\alpha$-content is upper semicontinuous, so passing to a subsequence if necessary,
    \begin{equation*}
        F=\lim_{n\to\infty}\bigl(r_n^{-1}(K-x)\cap B(0,1)\bigr)
    \end{equation*}
    satisfies $\mathcal{H}^\alpha_\infty(F)\geq 1$.
\end{proof}
Of course, we can combine the previous two results to obtain the following improvement of \cref{p:tangent-assouad}.
\begin{restatement}{ic:full-tan}
    Let $K$ be a compact set with $\dimA K=\alpha$.
    Then there is a weak tangent $F\in\Tan(K)$ with $\mathcal{H}^\alpha_\infty(F)\geq 1$.
\end{restatement}
\begin{proof}
    By \cref{p:tangent-assouad}, there is $E\in\Tan(K)$ such that $\mathcal{H}^\alpha(E)>0$.
    By \cite[Theorem~4.10]{zbl:0689.28003}, there is a compact $E'\subset E$ such that $0<\mathcal{H}^\alpha(E')<\infty$.
    Then by \cref{l:large-tan}, there is $F'\in\Tan(E')$ with $\mathcal{H}_\infty^\alpha(F')\geq 1$.
    But $F'\subset F$ for some $F\in\Tan(E)$, and by \cref{l:tan-compose}, $F\in\Tan(K)$ with $\mathcal{H}^\alpha_\infty(F)\geq\mathcal{H}^\alpha_\infty(F')\geq 1$.
\end{proof}

We now establish bounds on the pointwise Assouad dimension and tangents for general sets.
\begin{restatement}{it:gen-tangents}
    Let $K\subset\R^d$.
    Then:
    \begin{enumerate}[nl,r]
        \item\label{im:analy} If $K$ is analytic, for any $s$ such that $\mathcal{H}^s(K)>0$, there is a compact set $E\subset K$ with $\mathcal{H}^s(E)>0$ so that for each $x\in E$, there is a tangent $F\in\Tan(\overline{K},x)$ with $\mathcal{H}^{s}_\infty(F)\geq 1$.
        \item\label{im:compact} If $K$ is compact, there is an $x\in K$ such that $\dimA(K,x)\geq\dimuB K$.
    \end{enumerate}
\end{restatement}
\begin{proof}
    The proof of \cref{im:analy} follows directly from \cref{l:large-tan}, recalling that we can always find a compact subset $E\subset K$ such that $0<\mathcal{H}^s(E)<\infty$ (combine \cite[Theorem~8.19]{zbl:0819.28004} and \cite[Corollary~B.2.4]{zbl:1390.28012}).

    We now see \cref{im:compact}.
    Let $\dimuB K=t$.
    We first observe that for any $r>0$, there is an $x\in K$ so that $\dimuB B(x,r)\cap K=t$.
    In particular, we may inductively construct a nested sequence of balls $B(x_k,r_k)$ with $\lim_{k\to \infty}r_k=0$ so that $\dimuB K\cap B(x_k,r_k)=t$ for all $k\in\N$.
    Since $K$ is compact, take $x=\lim_{k\to\infty}x_k\in K$.
    We verify that $\dimA(K,x)\geq t$.
    Let $C>0$ and $\rho>0$ be arbitrary.
    Since the $x_k$ converge to $x$ and the $r_k$ converge to $0$, get some $k$ so that $B(x_k,r_k)\subset B(x,\rho)$.
    Thus for all $\epsilon>0$ and $r>0$ sufficiently small depending on $\epsilon$ and $\rho$, since $\dimuB K\cap B(x_k,r_k)=t$,
    \begin{equation*}
        N_r\bigl(B(x,\rho)\cap K\bigr)\geq N_r\bigl(B(x_k,r_k)\cap K\bigr)\geq C\left(\frac{r_k}{r}\right)^{t-\epsilon}.
    \end{equation*}
    Thus $\dimA(K,x)\geq t$.
\end{proof}

\begin{remark}
    Note that compactness is essential in \cref{it:gen-tangents}~\cref{im:compact} since there are sets with $\dimuB K=1$ but every point is isolated: consider, for instance, the set $E=\{(\log n)^{-1}:n=2,3,\ldots\}$.
    In this case, $\overline{E}=E\cup\{0\}$ and $\dimA(\overline{E},0)=1$.
    This example also shows that \cref{im:compact} can hold with exactly 1 point.
\end{remark}

Finally, we construct some general examples which go some way to showing that the results for general sets given in this section are sharp.
\begin{example}\label{ex:non-tan-attain}
    In general, the Assouad dimension can only be characterized by weak tangents rather than by tangents.
    For example, consider the set $K$ from \cite[Example~2.20]{zbl:1321.54059}, defined by
    \begin{equation*}
        K=\{0\}\cup\bigl\{2^{-k}+\ell 4^{-k}:k\in\N,\ell\in\{0,1,\ldots,k\}\bigr\}
    \end{equation*}
    Since $K$ contains arithmetic progressions of length $k$ for all $k\in\N$, $\dimA K=1$.
    However, $\dimA(K,x)=0$ for all $x\in K$ and, therefore, by \cref{p:pointwise-tangent}, $\dimH F = 0$ for all $F \in \Tan(K,x)$ and $x \in K$.
\end{example}

\begin{example}\label{ex:non-point-attain}
    We give an example of a compact set $K$ and a point $x\in K$ so that $\dimA(K,x)=1$ but each $F\in\Tan(K,x)$ consists of at most finitely many points.

    Set $a_k=4^{-k^2}$ and observe that $k a_{k+1}/a_k\leq 1/k$.
    For each $k\in\N$, write $\ell_k=\lfloor 2^k/k\rfloor$ and set
    \begin{equation*}
        K=\{0\}\cup\bigcup_{k=1}^\infty \left\{a_k\frac{2^k-\ell_k}{2^k},a_k\frac{2^k-\ell_k-1}{2^k},\ldots,a_k\right\}
    \end{equation*}
    and consider the point $x=0$.
    First observe for all $\epsilon>0$ and all $k$ sufficiently small depending on $\epsilon$,
    \begin{equation*}
        N_{2^{-k}\cdot a_k}\bigl(B(0,a_k)\cap K\bigr)\geq\frac{\ell_k}{2}\geq 2^{(1-\epsilon)k}
    \end{equation*}
    which gives that $\dimA(K,0)=1$.

    On the other hand, for $k\in\N$,
    \begin{equation*}
        a_k^{-1} K\cap B(0,1)\subset [0,a_{k+1}/a_k]\cup [1/k,1].
    \end{equation*}
    Since $k a_{k+1}/a_k\leq 1/k$, it follows that for any $\lambda\geq 1$ and $\lambda K\cap B(0,1)$ can be contained in a union of two intervals with arbitrarily small length as $\lambda$ diverges to $\infty$.
    Thus any tangent $F\in\Tan(K,0)$ consists of at most 2 points.
\end{example}

\subsection{Tangents of dynamically invariant sets} \label{eq:self-embeddable}
We recall from \cref{ex:non-tan-attain} and \cref{ex:non-point-attain} that the Assouad dimension of $K$ need not be attained as the Assouad dimension of a point, and even the Assouad dimension at a point need not be attained as the upper box dimension of a tangent at that point.

Now recall the definition of self-embeddability from \cref{d:self-embed}.
For self-embeddable sets, we can prove directly that the Assouad dimension of $K$ is attained as the Hausdorff dimension of a tangent.
In fact, the tangent can be chosen to have positive $\mathcal{H}^\alpha$-measure for $\alpha=\dimA K$.
\begin{theorem}\label{t:pointwise-tan}
    Let $K\subseteq\R^d$ be compact and self-embeddable with $\alpha=\dimA K$.
    Then there is a dense set of points $x\in K$ for which there exist $F\in\Tan(K,x)$ such that $\mathcal{H}^\alpha_\infty(F)\geq 2^{-\alpha}$.
    In particular,
    \begin{equation*}
        \dimP\{x\in K:\dimA(K,x)=\dimA K\}=\dimP K.
    \end{equation*}
\end{theorem}
\begin{proof}
    We first note that it suffices to construct a single point $x$ such that $\mathcal{H}^\alpha_\infty(F)\geq 2^{-\alpha}$.
    By self-embeddability and since $\dimA(K,x)=\dimA(f(K),f(x))$ for a bi-Lipschitz map $f$, this immediately yields a dense subset of such points.
    Moreover, recalling \cref{p:meas}~\cref{im:large-gdelta}, since $\dimA(K,x)\leq\dimA K$ for all $x\in K$, $\{x\in K:\dimA(K,x)=\dimA K\}$ is a dense $G_\delta$ subset of $K$ and therefore has packing dimension equal to the packing dimension of $K$ (see, for instance, \cite[Proposition~2.9]{zbl:1285.28011}).

    It therefore remains to construct such a point.
    Begin with an arbitrary ball $B(x_1,r_1)$ with $x_1\in K$ and $0<r_1\leq 1$.
    Since $K$ is self-embeddable, get a bi-Lipschitz map $f_1\colon K\to K\cap B(x_1,r_1)$.
    Since $\dimA f_1(K)=\alpha$, by \cref{ic:full-tan} there is a weak tangent $F_1$ of $f_1(K)$ such that $\mathcal{H}^\alpha_\infty(F_1)\geq 1$.
    Since $F_1$ is a weak tangent of $f_1(K)$, there is a similarity $T_1$ with similarity ratio $\lambda_1\geq 1$ such that $0\in T_1(K)$ and
    \begin{equation*}
        d_{\mathcal{H}}\bigl(T_1(f_1(K))\cap B(0,1),F_1\bigr)\leq 1.
    \end{equation*}
    Then choose $x_2\in K$ and $r_2\leq 1/2$ so that $B(x_2,r_2)\subset T_1^{-1}B^\circ(0,1)$.

    Repeating the above construction, next with the ball $B(x_2,r_2)$, and iterating, we obtain a sequence of similarity maps $(T_n)_{n=1}^\infty$ each with similarity ratio $\lambda_n\geq n$, bi-Lipschitz maps $f_n$, and compact sets $F_n$ such that
    \begin{enumerate}[nl]
        \item\label{im:incl} $T_{n+1}^{-1}B(0,1)\subseteq T_n^{-1}B(0,1)$,
        \item\label{im:conv} $\displaystyle d_{\mathcal{H}}\bigl(T_n(f_n(K))\cap B(0,1),F_n\bigr)\leq \frac{1}{n}$, and
        \item $\mathcal{H}^\alpha_\infty(F_n)\geq 1$.
    \end{enumerate}
    Let $x=\lim_{n\to\infty}T_n^{-1}(0)$ and note by \cref{im:incl} that $x\in T_n^{-1}B(0,1)$ for all $n\in\N$.
    Let $h_n$ be a similarity with similarity ratio $1/2$ such that
    \begin{equation*}
        d_{\mathcal{H}}\Bigl(\frac{\lambda_n}{2}(f_n(K)-x)\cap B(0,1),h_n(F_n)\Bigr)\leq\frac{1}{n}.
    \end{equation*}
    Observe that $\mathcal{H}^\alpha_\infty(h_n(F_n))\geq 2^{-\alpha}$.
    Thus passing to a subsequence if necessary, since $f_n(K)\subseteq K$, we may set
    \begin{equation*}
        F_0=\lim_{n\to\infty}\frac{\lambda_n}{2}(f_n(K)-x)\cap B(0,1)\qquad\text{and}\qquad F=\lim_{n\to\infty}\frac{\lambda_n}{2}(K-x)\cap B(0,1),
    \end{equation*}
    and observe that $F_0\subseteq F$.
    Again passing to a subsequence if necessary, by compactness of the orthogonal group, \cref{im:conv} and the triangle inequality, there is an isometry $h$ so that $\lim_{n\to\infty}h\circ h_n(F_n)=F_0$.
    Thus by upper semicontinuity of Hausdorff content,
    \begin{equation*}
        \mathcal{H}^\alpha_\infty(F)\geq\mathcal{H}^\alpha_\infty(F_0)\geq\lim_{n\to\infty}\mathcal{H}^\alpha_\infty(h\circ h_n(F_n))=2^{-\alpha}
    \end{equation*}
    as required.
\end{proof}
We recall from \cref{it:gen-tangents}~\cref{im:compact} that, for a general compact set $K$, the upper box dimension of $K$ provides a lower bound for the pointwise Assouad dimension at \emph{some} point.
For self-embeddable sets, we observe that the upper box dimension provides a uniform lower bound for the pointwise Assouad dimension at \emph{every} point.
On the other hand, the upper box dimension \emph{does not} lower bound the maximal dimension of a tangent.
\begin{proposition}\label{p:tan-lower}
    Let $K\subseteq\R^d$ be self-embeddable.
    Then for any $x\in K$, we have $\dimA(K,x)\geq\dimuB K$.
\end{proposition}
\begin{proof}
    Fix $\alpha<\dimuB K$ and $x\in K$.
    Let $C>0$ and $\rho>0$ be arbitrary.
    Since $K$ is self-embeddable, there is some bi-Lipschitz map $f\colon K\to B(x,\rho)$ so that $f(K)\subseteq K$.
    Since $\dimuB f(K)>\alpha$, there is some $0<r\leq\rho$ so that
    \begin{equation*}
        N_r(B(x,\rho)\cap K)\geq N_r(f(K))\geq C\Bigl(\frac{\rho}{r}\Bigr)^\alpha.
    \end{equation*}
    Since $C>0$ and $\rho>0$ were arbitrary, $\dimA(K,x)\geq\alpha$, as required.
\end{proof}

Now assuming uniform self-embeddability, we will see that the set of points with tangents that have positive $\mathcal{H}^{\alpha}$-measure has full Hausdorff dimension for $\alpha=\dimA K$.
Since uniformly self-embeddable sets satisfy the hypotheses of \cite[Theorem~4]{zbl:0683.58034}, it always holds that $\dimuB K=\dimH K$ (see also \cite[Theorem~2.10]{zbl:1305.28021}).
On the other hand, it can happen in this class of sets that $\dimuB K<\alpha$: for example, this is the situation for self-similar sets in $\R$ with $\dimuB K<1$ which fail the weak separation condition; see \cite[Theorem~1.3]{zbl:1317.28014}.
We provide a subset of full Hausdorff dimension for which each point has a tangent with positive Hausdorff $\alpha$-measure.

The idea of the proof is essentially as follows.
Let $F$ be a weak tangent for $K$ with strictly positive Hausdorff $\alpha$-content.
For each $s<\dimuB K$, using the implicit method of \cite[Theorem~4]{zbl:0683.58034}, we can construct a well-distributed set of $N$ balls at resolution $\delta$, where $\delta^{-s}\ll N$.
Then, inside each ball, using uniform self-embeddability, we can map an image of an approximate tangent $T_\delta^{-1}(B(0,1))\cap K\approx F$ where $T_\delta$ has similarity ratio $\lambda$.
Choosing $N$ to be large, the resulting collection of images of the approximate tangent $F$ is again a family of well-distributed balls at resolution $\lambda^{-1}\delta$, with $(\lambda^{-1}\delta)^{-s}\approx N$.
Repeating this construction along a sequence of tangents converging to $F$ yields a set $E$ with $\dimH E\geq s$ such that each $x\in E$ has a tangent which is an image of $F$ (up to some negligible distortion), which has positive Hausdorff $\alpha$-content by upper semicontinuity of content.

We fix a compact set $K$.
To simplify notation, we say that a function $f\colon K\to K$ is in $\mathcal{G}(z,r,c)$ for $z\in K$ and $c,r>0$ if $f(K)\subset B(z,r)$ and
\begin{equation*}
    cr|x-y|\leq|f(x)-f(y)|\leq c^{-1}r|x-y|
\end{equation*}
for all $x,y\in K$.
\begin{theorem}\label{t:uniformly-self-embeddable}
    Let $K\subset\R^d$ be uniformly self-embeddable  and let $\alpha=\dimA K$.
    Then
    \begin{equation*}
        \dimH\{x\in K:\exists F\in\Tan(K,x)\text{ with }\mathcal{H}^\alpha_\infty(F)\gtrsim 1\}=\dimH K = \dimuB K.
    \end{equation*}
\end{theorem}
\begin{proof}
    Write $\alpha=\dimA K$.
    If $\dimuB K=0$ we are done; otherwise, let $0<s<\dimuB K$ be arbitrary.
    Since $K$ is uniformly self-embeddable, there is a constant $a\in(0,1)$ so that for each $z\in K$ and $0<r\leq\diam K$ there is a map $f\in\mathcal{G}(z,r,a)$.
    Next, from \cref{ic:full-tan}, there is a compact set $F\subset B(0,1)$ with $\mathcal{H}^\alpha_\infty(F)\geq 1$ and a sequence of similarities $(T_k)_{k=1}^\infty$ with similarity ratios $(\lambda_k)_{k=1}^\infty$ such that
    \begin{equation*}
        F=\lim_{k\to\infty}T_k(K)\cap B(0,1)
    \end{equation*}
    with respect to the Hausdorff metric.
    Set $Q_k=T_k^{-1}(B(0,1))\cap K$.
    We will construct a Cantor set $E\subset K$ of points each of which has pointwise Assouad dimension at least $\alpha$ and has $\dimH E\geq s$.

    We begin with a preliminary construction.
    First, since $s<\dimuB K$, there is some $r_0>0$ and a collection of points $\{y_i\}_{i=1}^{N_0}\subset K$ such that $|y_i-y_j|>3 r_0$ for all $i\neq j$ and $N_0\geq 2^sa^{-s}r_0^{-s}$.
    Now for each $i$, take a map $\phi_i\in\mathcal{G}(y_i,r_0,a)$.
    Write $\mathcal{I}=\{1,\ldots,N_0\}$, and for $\mtt{i}=(i_1,\ldots,i_n)\in\mathcal{I}^n$ set
    \begin{equation*}
        \phi_{\mtt{i}}=\phi_{i_1}\circ\cdots\circ\phi_{i_n},
    \end{equation*}
    and, having fixed some $x_0\in K$, write $x_\mtt{i}=\phi_{\mtt{i}}(x_0)\in\phi_{\mtt{i}}(K)$.
    Observe that if the maximal length of a common prefix of $\mtt{i}$ and $\mtt{j}$ is $m$, then
    \begin{equation*}
        \dist(\phi_{\mtt{i}}(K),\phi_{\mtt{j}}(K))\geq r_0 (ar_0)^m.
    \end{equation*}

    We now begin our inductive construction.
    Without loss of generality, we may assume that $\lambda_n\geq 12$ for all $n\in\N$ and $r_0\leq 1$.
    First, for each $n\in\N$, define constants $(m_n)_{n=1}^\infty\subset\{0\}\cup\N$ and $(\rho_n)_{n=1}^\infty$ converging monotonically to zero from above by the rules
    \begin{enumerate}[nl]
        \item $\displaystyle 2^{-m_{n}}\leq\frac{a^{2}r_0\lambda_n^{-1}}{3}$,
        \item $\rho_0=\diam K$, and
        \item $\displaystyle\rho_{n}\coloneqq \rho_{n-1}\cdot\frac{a\lambda_n^{-1}\cdot(ar_0)^{m_{n}}}{3}$.
    \end{enumerate}
    Next, for $n\in\N\cup\{0\}$ we inductively choose points $y_{n,\mtt{i}}\in K$ and maps $\Psi_{n,\mtt{i}}\in\mathcal{G}(y_{n,\mtt{i}},\rho_n,a)$ for $\mtt{i}\in\mathcal{I}^{m_1}\times\cdots\times\mathcal{I}^{m_n}$.
    Let $\varnothing$ denote the empty word and let $y_{0,\varnothing}\in K$ be arbitrary and let $\Psi_{0,\varnothing}$ denote the identity map.
    Then for each $\mtt{k}=\mtt{i}\mtt{j}$ with $\mtt{i}\in\mathcal{I}^{m_1}\times\cdots\times\mathcal{I}^{m_{n-1}}$ and $\mtt{j}\in\mathcal{I}^{m_{n}}$, sequentially choose:
    \begin{enumerate}[nl, resume]
        \item $\psi_{n,\mtt{k}}\in\mathcal{G}(\Psi_{n-1,\mtt{i}}(x_\mtt{j}),\rho_{n}\lambda_{n} a^{-1},a)$
        \item $y_{n,\mtt{k}}= \psi_{n,\mtt{k}}\circ T_n^{-1}(0)$
        \item $\Psi_{n,\mtt{k}}\in\mathcal{G}(y_{n,\mtt{k}},\rho_n,a)$
    \end{enumerate}
    Finally, write $\mathcal{J}_0=\{\varnothing\}$, $\mathcal{J}_n=\mathcal{I}^{m_1}\times\cdots\times\mathcal{I}^{m_n}$ for $n\in\N$, and let
    \begin{equation*}
        E_n=\bigcup_{\mtt{i}\in\mathcal{J}_n}B(y_{n,\mtt{i}\mtt{j}},3\rho_n)\qquad\text{and}\qquad E=K\cap\bigcap_{n=1}^\infty E_n.
    \end{equation*}
    Suppose $\mtt{i}\in\mathcal{J}_{n-1}$ and $\mtt{j}\in\mathcal{I}^{m_n}$.
    Since $x_\mtt{j}\in K$, $\Psi_{n-1,\mtt{i}}(K)\subset B(y_{n-1,\mtt{i}},\rho_{n-1})$, and $y_{n,\mtt{i}\mtt{j}}\in\psi_{n,\mtt{i}\mtt{j}}(K)\subset B(\Psi_{n-1,\mtt{i}}(x_\mtt{j}),\rho_{n})$, we conclude since $\rho_{n-1}\geq 3\rho_n$ that
    \begin{equation*}
        B(y_{n,\mtt{i}\mtt{j}},3\rho_n)\subset B(y_{n-1,\mtt{i}},3\rho_{n-1}).
    \end{equation*}
    Moreover, $y_{n,\mtt{i}\mtt{j}}\in K$, so the sets $E_n$ are non-empty nested compact sets and therefore $E$ is non-empty.

    We next observe the following fundamental separation properties of the balls in the construction of the sets $E_n$.
    Let $n\in\N$ and suppose $\mtt{j}_1\neq\mtt{j}_2$ in $\mathcal{I}^{m_n}$ and $\mtt{i}\in\mathcal{J}_{n-1}$ (writing $\mathcal{J}_0=\{\varnothing\}$).
    Suppose $\mtt{j}_1$ and $\mtt{j}_2$ have a common prefix of maximal length $m$.
    First recall that $|x_{\mtt{j}_1}-x_{\mtt{j}_2}|\geq r_0(ar_0)^{m}$, so that
    \begin{equation*}
        |\Psi_{n-1,\mtt{i}}(x_{\mtt{j}_1})-\Psi_{n-1,\mtt{i}}(x_{\mtt{j}_2})|\geq \rho_{n-1}(ar_0)^{m+1}.
    \end{equation*}
    Then, since for $j=1,2$
    \begin{equation*}
        y_{n,\mtt{i}\mtt{j}_j}\in\psi_{n,\mtt{i}\mtt{j}_j}(K)\subset B\Bigl(\Psi_{n-1,\mtt{i}}(x_{\mtt{j}_j}),\frac{\rho_{n-1}(ar_0)^{m_{n}}}{3}\Bigr)
    \end{equation*}
    we observe that
    \begin{equation*}
        |y_{n,\mtt{i}\mtt{j}_1}-y_{n,\mtt{i}\mtt{j}_2}|\geq \rho_{n-1}(ar_0)^{m+1}-2\frac{\rho_{n-1}(ar_0)^{m_{n}}}{3}\geq\frac{\rho_{n-1}(ar_0)^{m+1}}{3}.
    \end{equation*}
    Since we assumed that $\lambda_n\geq 12$, by the triangle inequality
    \begin{equation}\label{e:ball-dist}
        \dist\bigl(B(y_{n,\mtt{i}\mtt{j}_1},3\rho_n),B(y_{n,\mtt{i}\mtt{j}_2},3\rho_n)\bigr)\geq \frac{\rho_{n-1}(ar_0)^{m+1}}{3}-6\rho_n\geq\frac{\rho_{n-1}(ar_0)^{m+1}}{6}.
    \end{equation}

    We first show that $\dimH E\geq s$.
    By the method of repeated subdivision, define a Borel probability measure $\mu$ with $\supp\mu=E$ and for $\mtt{i}\in\mathcal{J}_n$,
    \begin{equation*}
        \mu(B(y_{n,\mtt{i}},3\rho_n)\cap K)=\frac{1}{\#\mathcal{J}_n}.
    \end{equation*}
    Now suppose $U$ is an arbitrary open set with $U\cap E\neq\varnothing$.
    Intending to use the mass distribution principle, we estimate $\mu(U)$.
    Assuming that $U$ has sufficiently small diameter, let $n\in\N$ be maximal so that
    \begin{equation*}
        \diam U\leq \frac{a^{-1}\lambda_n}{2}\rho_n=\frac{\rho_{n-1}(ar_0)^{m_n}}{6}.
    \end{equation*}
    By \cref{e:ball-dist}, there is a unique $\mtt{i}\in\mathcal{J}_{n}$ such that $U\cap B(y_{n,\mtt{i}},3\rho_n)\neq\varnothing$.
    We first recall by choice of the constants $m_n$ that
    \begin{align*}
        \rho_n&=(\diam K)\cdot\Bigl(\frac{a^2 r_0}{3}\Bigr)^{n}\lambda_1^{-1}\cdots \lambda_n^{-1}(ar_0)^{m_1+\cdots+m_n}\\
              &\geq(\diam K) 2^{-(m_1+\cdots+m_n)}(ar_0)^{m_1+\cdots+m_n}.
    \end{align*}
    There are two cases.
    First assume $\rho_n/6<\diam U$.
    Thus
    \begin{align*}
        \mu(U)&\leq\frac{1}{\#\mathcal{J}_n}\leq\Bigl(\frac{1}{2}ar_0\Bigr)^{s(m_1+\cdots+m_n)}\leq(\diam K)^{-s}\rho_n^s\\
              &\leq \Bigl(\frac{6}{\diam K}\Bigr)^s\cdot(\diam U)^s.
    \end{align*}
    Otherwise, let $k\in\{0,\ldots,m_{n+1}-1\}$ be so that
    \begin{equation*}
        \frac{\rho_n (ar_0)^{k+1}}{6}<\diam U\leq\frac{\rho_n (ar_0)^{k}}{6}.
    \end{equation*}
    By \cref{e:ball-dist}, $U$ intersects at most $N_0^{m_n-k}$ balls $B(y_{n+1,\omega},3\rho_{n+1})$ for $\omega\in\mathcal{J}_{n+1}$, so since $2^{-sk}\leq 1$,
    \begin{align*}
        \mu(U)&\leq\frac{1}{\#\mathcal{J}_n\cdot N_0^k}\leq(\diam K)^{-s}\rho_n^s\cdot(2^{-s}(ar_0)^s)^k\\
              &\leq\Bigl(\frac{6}{ar_0\diam K}\Bigr)^s\cdot\Bigl(\frac{\rho_n(ar_0)^{k+1}}{6}\Bigr)^s\\
              &\leq\Bigl(\frac{6}{ar_0\diam K}\Bigr)^s\cdot\bigl(\diam U\bigr)^s.
    \end{align*}
    This treats all possible small values of $\diam U$, so there is a constant $M>0$ such that $\mu(U)\leq M(\diam U)^s$.
    Thus $\dimH E\geq s$ by the mass distribution principle.

    Now fix
    \begin{equation*}
        C=(3+a^{-2})^{-\alpha}.
    \end{equation*}
    We will show that each $z\in E$ has a tangent with Hausdorff $\alpha$-content at least $C$.
    Let $z\in E$ and define
    \begin{equation*}
        S_n(x)=\frac{x-z}{\rho_n(3+a^{-2})}.
    \end{equation*}
    Our tangent will be an accumulation point of the sequence $(S_n(K)\cap B(0,1))_{n=1}^\infty$.
    Now fix $n\in\N$.
    Since $z\in E$, there is some $\omega\in\mathcal{J}_n$ so that $z\in B(y_{n,\omega},3\rho_n)$.
    By choice of $y_{n,\omega}$, $Q_n=B\bigl(\psi_{n,\omega}^{-1}(y_{n,\omega}),\lambda_n^{-1}\bigr)\cap K$ so that
    \begin{equation*}
        \psi_{n,\omega}(Q_n)\subseteq B\bigl(y_{n,\omega},\rho_n a^{-2}\bigr)\cap K\subseteq B\bigl(z,\rho_n(3+a^{-2})\bigr)\cap K
    \end{equation*}
    and therefore, writing $\Phi_n=S_n\circ\psi_{n,\omega}\circ T_n^{-1}$,
    \begin{equation*}
        \Phi_n(T_n(K)\cap B(0,1))\subset S_n(K)\cap B(0,1).
    \end{equation*}
    Then for $x,y\in T_n(K)\cap B(0,1)$, by the choice of $\psi$ in (4),
    \begin{equation}\label{e:phi-dist-bounds}
        \frac{|x-y|}{3+a^{-2}}\leq|\Phi_n(x)-\Phi_n(y)|\leq\frac{|x-y|}{a^2(3+a^{-2})}.
    \end{equation}
    Now, passing to a subsequence $(n_k)_{k=1}^\infty$, we can ensure that
    \begin{equation*}
        \lim_{k\to\infty}\Phi_{n_k}(F)=Z_0\qquad\text{and}\qquad\lim_{k\to\infty}S_{n_k}(K)\cap B(0,1)=Z.
    \end{equation*}
    Moreover, recall that $\lim_{k\to\infty}T_{n_k}(K)\cap B(0,1)=F$ and $\mathcal{H}^\alpha_\infty(F)\geq 1$.
    Observe by \cref{e:phi-dist-bounds} that $\mathcal{H}^\alpha_\infty(\Phi_{n_k}(F))\geq C$ for each $k$, so by upper semicontinuity of Hausdorff content, $\mathcal{H}^\alpha_\infty(Z_0)\geq C$.
    But again by \cref{e:phi-dist-bounds},
    \begin{align*}
        d_{\mathcal{H}}\bigl(Z_0,\Phi_{n_k}(T_{n_k}(K)\cap B(0,1))\bigr)\leq d_{\mathcal{H}}(Z_0,\Phi_{n_k}(F))+\frac{d_{\mathcal{H}}\bigl(F,T_{n_k}(K)\cap B(0,1)\bigr)}{a^2(3+a^{-2})}
    \end{align*}
    so in fact $Z_0\subset Z$ and $\mathcal{H}^\alpha_\infty(Z)\geq C$, as claimed.
\end{proof}
\begin{remark}
    We note that the upper distortion bound in the definition of uniform self-embeddability is used only at the very last step to guarantee that the images $\Phi_{n_k}(T_{n_k}(K)\cap B(0,1))$ converge to a large set whenever the $T_{n_k}(K)\cap B(0,1)$ converge to a large set.
\end{remark}

\begin{acknowledgements}
    The authors thank Roope Anttila for interesting discussions around the topics in this document.
    They also thank Amlan Banaji for some comments on a draft version of this document.
    Finally, they thank Balázs Bárány and Lars Olsen for many valuable comments on this document, and in particular pointing out that \cref{e:pack-full} follows from density.
    AR was supported by EPSRC Grant EP/V520123/1 and the Natural Sciences and Engineering Research Council of Canada.
    This project began while AR visited AK at the University of Oulu, which was funded by a scholarship from the London Math Society and the Cecil King Foundation.
    He thanks the various members of the department for their hospitality during the visit.
\end{acknowledgements}
\end{document}